\newtheorem{theorem}{Theorem}

\newtheorem{proposition}[theorem]{Proposition}

\newtheorem{corollary}[theorem]{Corollary}

\documentclass{article}
\usepackage{amsfonts} 
\usepackage{amsmath} 
\usepackage[subnum]{cases}
\usepackage{ifpdf}

\title{Monotonicity of a Class of Integral Functionals}

\author{Stefano Bertoni\footnote{Dipartimento di Matematica, Universit\`a di Trento (Italy). E-mail: bertoni@science.unitn.it} 
}

\begin{document}
\maketitle 
\begin{abstract}
In this note we prove a condition of monotonicity for the integral functional 
$ F(g) = \int_a^b h(x)\, d[-g(x)] $ with respect to $g$, a function of bounded variation. 
\end{abstract}



\section{Introduction}
In the article  \cite{bertoni:eq} (``\emph{Nontrivial Equilibria of a Quasilinear Population Model}'', in progress), I study a functional $R(u)$ ($u\in L^1(0,\infty)$), said \emph{generalized net reproduction rate},  
to prove existence of non--zero equilibria in a general structured population model.

The monotonicity of $R(u)$ is used in a Corollary to prove the non-existence of a non--zero stationary population if $R(0) <1$ (a sufficient condition of existence being $R(0)>1$). 

The original proposition about monotonicity, not so immediate, will be reduced to the integration by parts of an improper Stieltjes integral: 
$$
\int_a^\infty h(x)\, d[-g(x)]  = h(a)\,g(a) - \lim_{b\to\infty} h(b)\, g(b) + \int_a^\infty g(x)\, dh(x)
$$

\section{Monotonicity Propositions}
Assume  $0< a < b \leq \infty$. 

From now on we denote via $G(b) $ the value of $G(b) $ if $b < \infty$ and 
$\lim_{x\to\infty} G(x) $ if $b = \infty$. 
I will denote respectively in the cases $ [a, b] $ and $ [a, \infty) $.  

\begin{proposition}
\label{Prop1}
Let $H$, $G$ be two given functions on $I$. 

Let $H$ be increasing (non-decreasing), bounded, non-negative. 
Let $G$ be continuous and of bounded variation.  

Define 
\begin{equation}
{\cal F}(G) := \int_a^b H(x)\,d[-G(x)] . 
\end{equation}
If $G(b) = 0$, then ${\cal F}$ is increasing (non-decreasing) with respect to $G$, i.e.  
let be $A := \{\phi |\phi \in C([a,b]) \cap BV[a,b], \phi(b) = 0\}$:  
if $ G_1, G_2 \in A $ and  $ G_1 < G_2$ , then $ {\cal F}(G_1) < {\cal F}(G_2)$ 
(respectively $ {\cal F}(G_1) \leq {\cal F}(G_2)$). 
\end{proposition}
\paragraph{\it Proof.} 
\emph{a)} Consider first the case $b < \infty$. ${\cal F}(G) $  is well--defined; integrating by parts we have: 
\begin{equation}
{\cal F}(G) = -H(b)\,G(b) + H(a)\,G(a) + \int_a^b G(x)\,dH(x) 
= H(a)\,G(a) + \int_a^b G(x)\,dH(x) .
\end{equation}
The conclusion is immediate. 
\par\medskip\noindent 
\emph{b)} Consider the case $b=\infty$. For $H$ bounded and $G(x) $ converging for $x \to \infty$ we obtain immediately the existence of the improper integral and extend the formula of case a).  

If  $H(x) $ is not strictly increasing but only non--decreasing, the functional $\cal F$ is only non--decreasing with respect to $G$. 

\begin{corollary}\label{cor1}
Let   $H$, $G$ given functions on $I$. 

Let $H$ be decreasing (non--increasing), bounded, non-negative. 
Let $G $ be continuous and of bounded variation.

Define ${\cal F}_0 (G) := \int_a^b H(x)\,  dG(x)   $.

If $G(b) = 0 $, 
then $ {\cal F}_0$ is increasing (non--decreasing) with respect to $G$. 
\end{corollary}
\paragraph{\it Example 1.} Consider the functional 
\begin{equation}
{\cal I}(f) = \int_0^\infty dx\, h(x)\, f(x)\, e^{-\int_0^x dy\, f(y)}
\end{equation}
where $h$ is positive, increasing and bounded. If  $ f \in L^1_{\it loc}(0, \infty) $, $ f\geq 0 $
and  
$\int_0^\infty dy f(y) = \infty$ 
($f\not\in L^1(0,\infty) $),  
then  $\cal I$ is decreasing with respect to $f$. 

This is a particular case of Prop. \ref{Prop1}, where $g(x) = e^{-\int_0^x dy\, f(y)} $
and $${\cal I}(f) = \int_0^\infty dx \, h(x)\, d[- e^{-\int_0^x dy\, f(y)} ] . $$

\begin{corollary}\label{cor2} 
Consider  $ u\in L^1 (0,\infty) $ and the functional 
\begin{equation}
 R(u) = \int_0^\infty h(x, u(\cdot))\, f(x, u(\cdot))\, e^{-\int_0^x dy\, f(y, u(\cdot))} 
\end{equation}
where $h$  and $f$ are  defined from $ (0,\infty) \times L^1_+(0,\infty) $ in $ [0,\infty) $, $h$ is positive and bounded, $x\mapsto f \in L^1_{loc}(0, \infty) $ and $\int_0^\infty dy\,f(y) = 0$, and 
\begin{itemize}
\item
 let $x \mapsto h(x, u) $ be non-decreasing (increasing) for fixed  $ u $ 
\item 
$ u \mapsto h(x, u) $ decreasing (o non--increasing) for fixed $ x  $ 
\item 
$ u \mapsto f(x,u) $ non-decreasing (o increasing) for fixed $ x $  
\end{itemize}
Then $ R(u) $ is decreasing with respect to u. 
\end{corollary}
Proof. Take  $u_1, u_2 \in L^1_+(0,\infty) $ with  $u_1 < u_2$. 
For Proposition 1, the integral 
$$
\int_0^\infty h(x, u_1)\, f(x, u)\,  e^{-\int_0^x dy	\, f(y, u)} 
$$
is decreasing with respect to $ f $, that is non--decreasing in $u$: therefore this integral is non--increasing in $ u $  and we have  
\begin{eqnarray}
\int_0^\infty h(x, u_1)\, f(x, u_1)  e^{-\int_0^x dy\, f(y, u_1)}  
\geq 
\int_0^\infty h(x, u_1)\,f(x, u_2)  e^{-\int_0^x dy\, f(y, u_2)}.
\end{eqnarray}
As $ f $ is decreasing with respect to $ u$, we have  
\begin{eqnarray}
\int_0^\infty h(x, u_1)\, f(x, u_2) \, e^{-\int_0^x dy\, f(y, u_2)} >
 \int_0^\infty h(x, u_2)\, f(x, u_2) \, e^{-\int_0^x dy\, f(y, u_2)}, 
\end{eqnarray}

so that  $R(u_1) > R(u_2) $. 

(The case of the alternative conditions, given by the parenthesis, is analogous). 

\paragraph{\it Example 2.} Corollary \ref{cor2} is applied to a model of population dynamics: let $ u = u(t, x) \geq 0$ be a population density with respect to age or size $x \geq 0$. Existence of stationary solutions (i. e. equilibria) $ u=u(x) $  is related to a functional $ R(u) $, the net reproduction rate. In a generalized model (see \cite{bertoni:eq}) where $g$ and $\mu $  depend on $u$ in an infinite--dimensional kind,  $R(u)$ is represented by  
\begin{equation}
R(u) = \int_0^\infty dx\, \beta(x, u(\cdot))\, \frac{e^{-\int_0^x dy\, \frac{\mu(y, u(\cdot))}{g(y,u(\cdot))}}}{g(x,u(\cdot))}
\end{equation}
where  $  \beta $  represents fertility, $ \mu $   mortality and  $  g  $  is a coefficient of growth (the detailed model is given and discussed in \cite{bertoni:eq}). 

The condition of existence of a nonzero steady solution (with suitable regularity conditions) is requiring that $  R(u) = 1  $; see \cite{casa:amod,casa:bas} and \cite{bertoni:eq}. See also \cite{casn:sta,faha:sta,nk:agm}.

If  $  R(0) < 1 $   and monotonicity conditions hold, the zero solution is the unique equilibrium. 

I prove in  \cite{bertoni:eq} that  $  R(0)>1  $  is a sufficient condition for existence of nontrivial stationary solutions. If monotonicity conditions do not hold, then  $  R(0)>1  $  is sufficient but it is not necessary and it is simple to give a counterexample. 

\section{More about the Application} 
The model is a generalized version of the classic Lotka-MacKendrick population model: consider a population density $u=u(t,x)$, where  $t \in [0,T]$ represents \emph{time}, $x\in (0,\infty)$ is \emph{age} or \emph{size} and the total population $P(t)$ is  
$$ P(t) = \int_0^\infty\, u(t,x)\, dx. $$ 

Consider the following functions: growth/diffusion  $g= g(x, u) $,  mortality  $\mu=\mu(t,u)$, fertility $\beta=\beta(x, u)$, depending on $x$ and infinite--dimensionally depending on the population density $u(t, \cdot)$. The model is 
\begin{eqnarray}
\label{diff} && 
	u_t(t,x) + ( g(x, u(t, \cdot))\,u(t,x) )_x + \mu(x, u(t, \cdot) )\,u(t,x) = 0 , 
	\\  && 
	\label{newborns}
	g(0,u(t, \cdot))\,u(t,0) = \int_0^\infty dx\, \beta(x,u(t, \cdot))\,u(t, x) . 
\end{eqnarray} 
In particular, Eq. (\ref{newborns}) gives the newborns. 

The \emph{generalized net reproduction rate} is defined as 
\begin{equation}
	R(u) = \int_0^\infty \beta(x,u)\, \Pi(x,u) \,dx ,
\end{equation} 
where 
$\displaystyle 
\Pi(x,u) =\frac{1}{g(x,u)}\,e^{-\int_0^x \frac{\mu(y, u)}{g(y, u)} dy} $ is an auxiliary function, said \emph{generalized survival probability} and it represents a stationary solution of Eq. (\ref{diff}), i. e. the differential part of the model. 

In general $\beta$ and $ \Pi$ depend on $u$ in a functional way: for instance in Calsina and Saldana \cite{casa:amod,casa:bas} the dependence is given through a weighted integral; in my paper \cite{bertoni:eq} the dependence is  infinite-dimensional in a more general way, to manage hierarchical models. 

Some examples are populations where fertility or mortality are influenced only by the immediately superior size: for instance a population of trees in a forest, where the contended resource is the light, that is intercepted by immediately taller trees than trees of size $x$ but not by the trees that are very taller than $x$. (For a case of tree population model, see \cite{kraev}).

A stationary solution $u$ of (\ref{diff})--(\ref{newborns}) exists if and only if $u$ 
satisfies the functional equation 
\begin{equation}
\label{stat-eq}
	u = G(u)\,\Pi(u), 
\end{equation}
where $G(u(\cdot))=\int_0^\infty \beta(x', u(\cdot))\,u(x')\,dx'$. 

Eq. (\ref{stat-eq}) is related to the condition $R(u) = 1$ that is used to prove the existence of nontrivial stationary solution (that is, nonzero). 
Under suitable regularity conditions, we have that $R(0)>1$ is a sufficient condition. 

With additional conditions on monotonicity of $\beta/g$ and $\mu/g$, the reproduction rate $R(u)$ is monotone decreasing and we exclude existence of nontrivial solution if $R(0) < 1$. This is is a recurrent condition in dynamics of populations.

\section{Other Recurrences of the Functional in Literature}

Conditions on $H$ and $G$ in  Prop. 1 are analogous to conditions given in \cite{HeMa:wim}, Teorema 2.1, b) 
Teorema \cite{HeMa:wim} Let $ -\infty < a  < b \leq \infty$ and let $h$ and $g$ be positive functions on $ (a, b) $, where $g$ is continuous on  $ (a, b) $. 

Assume that $h$ is increasing on  
$(a, b) $ and $g$ is decreasing on $ (a, b) $ where  $ g (b^-) = 0$. Then,  for any $p\in (0, 1] $, 
\begin{equation}
\int_a^b h(x) d[-g(x)] \leq  \left( \int_a^b h^p(x) d[-g^p(x)] \right)^p \qquad\qquad {\rm(1.2)}
\end{equation}
If $1\leq p<\infty$, then the inequality (1.2) holds in the reversed direction.

In \cite{PePePe:iimf}, the theorem above extends from $t^p$ to concave and convex functions $\phi$, when they are positive and differentiable.  

At the present I have no ideas if this fact would have any meaning for $R(u) $ or eventually estimates of it in the spaces $L^p$, however I think that the similarities of conditions is not a coincidence. 

Heinig and Maligranda's original paper \cite{HeMa:wim}  treats monotone functions and H\"older inequalities on Hardy spaces. 
A related field can be about Fredholm-Volterra equations.

\end{document}